\documentclass[12pt]{article}

\usepackage{latexsym,amssymb,amsmath}
\usepackage[dvips]{graphicx}
\newcommand{\C}{\mathbb C}
\newcommand{\Z}{\mathbb Z}
\newcommand{\N}{\mathbb N}
\newcommand{\R}{\mathbb R}

\newcommand{\slt}{SL(2,\C)}

\newtheorem{lem}{Lemma}[section]

\newtheorem{co}[lem]{Corollary}
\newtheorem{thm}[lem]{Theorem}
\newtheorem{prop}[lem]{Proposition}

\newenvironment{proof}{\textbf{Proof.}}{\newline\hspace*{\fill}{$\Box$}\\}

\begin{document}
\title{Strictly ascending HNN extensions of finite rank free groups that
are linear over $\Z$}
\author{J.\,O.\,Button\\
Selwyn College\\
University of Cambridge\\
Cambridge CB3 9DQ\\
U.K.\\
\texttt{jb128@dpmms.cam.ac.uk}}
\date{}
\maketitle
\begin{abstract}
We find strictly ascending HNN extensions of finite rank free groups
possessing a presentation 2-complex which is a non positively curved
square complex. On showing these groups are word hyperbolic, we have
by results of Wise and Agol that they are linear over $\Z$. An example
is the endomorphism of the free group on $a,b$ that sends $a$ to
$ab^{-1}a^2b$ and $b$ to $ba^{-1}b^2a$.
\end{abstract}

\section{Introduction}

The work of Agol \cite{ag} and Wise \cite{wslng} was able to solve the
remaining questions on fundamental groups of closed hyperbolic 3-manifolds,
and even to provide new insights such as these groups are all linear
over $\Z$, but no 3-manifold theory is invoked in these proofs. Instead
the only hypotheses required are that we have a word hyperbolic group
which acts properly and cocompactly on a CAT(0) cube complex. Here a
cube complex is formed by gluing together cubes of various dimensions
and a purely combinatorial criterion can be used to determine if such
a complex admits a metric of non positive curvature. A CAT(0) cube
complex is then a simply connected cube complex admitting such a metric.
Therefore if we are able to satisfy these two conditions for other classes
of groups then the conclusions will hold here too. Alternatively we might
try and find a class of word hyperbolic groups where one or more of the
conclusions do not hold, and this would show that such groups cannot be
``cubulated''.

In this paper we look at groups that in several senses lie near to
closed hyperbolic 3-manifold groups. One such class is that of
$F_n$ by $\Z$ groups, where $F_n$ denotes the free group of finite rank $n$.
Such groups are obtained by taking an automorphism $\alpha$ of $F_n$ and
forming the semidirect product $F_n\rtimes_\alpha\Z$, in analogy with groups
$S_g\rtimes_\alpha\Z$ for $S_g$ the fundamental group of a closed orientable
surface of genus $g$. We know from the work of Thurston that if one of
these latter groups 
is word hyperbolic then it will be the fundamental group of a closed
hyperbolic 3-manifold. Therefore we could look at word hyperbolic 
$F_n$ by $\Z$ groups and see if we can always find a CAT(0) cube complex
on which such a group acts properly and cocompactly. However this is a
major result for closed hyperbolic 3-manifold groups and requires use of the
Kahn-Markovic theorem \cite{km} on the existence of many surface subgroups,
as well as \cite{bws} which uses these codimension 1 subgroups to form the 
CAT(0) complex on which the group acts.

In this paper we adopt a far more direct and basic approach which says that
if we are given a finitely presented group $G$ then the simplest CAT(0)
cube complex we could hope for where $G$ acts properly and cocompactly
would be a 2 dimensional complex made only out of squares. Indeed the
non positively curved condition takes on a particularly concise form here
which is easily checked. As for the action, this would immediately follow
if we could find a square complex which is a presentation 2-complex for $G$,
as then we can take the universal cover. Thus we need to find a presentation
for $G$ where every relator has length 4 (giving us the square complex), then
look at the various subwords of length 2 appearing in these relators, inverses
and cyclic permutations (which gives rise to the non positively curved 
condition). Such presentations have appeared before, for instance knot groups,
but this process is only of use to us if our presentation gives rise to
a word hyperbolic group, which is needed for the Wise - Agol machine to be
applied.

Here our class of groups is closely related to $F_n$ by $\Z$ groups but rather
than forming an HNN extension using the automorphism $\alpha$, we use an
injective (but not surjective) endomorphism $\theta$ and we refer to these
HNN extensions as strictly ascending. These groups have a lot in common
with $F_n$ by $\Z$ groups and can be word hyperbolic. However in the word
hyperbolic case there appear to be no examples in the literature of
strictly ascending HNN extensions of finite rank free groups (or even just
of finitely generated groups) which have been proved to be linear over $\C$,
let alone over $\Z$. Now linearity over $\Z$ is a very strong property and 
the nested sequence of conjugates of the base group occurring in a strictly
ascending HNN extension of a finitely generated group suggests that such a
group cannot be linear over $\Z$ (for instance the group
$\langle a,t|tat^{-1}=a^2\rangle$ with base $\Z\cong F_1$ is linear over
$\C$ but not $\Z$). However here we give specific examples of finite
presentations defining strictly ascending HNN extensions of finite rank free
groups which have another presentation where the universal cover of this
presentation 2-complex is a CAT(0) cube complex. 
We then establish word hyperbolicity,
allowing application of Agol - Wise to obtain strictly ascending HNN extensions
which are linear over $\Z$.

In Section 2 we put together some standard material which allows us to show
that $F_n$ by $\Z$ groups and strictly ascending HNN extensions of finite
rank free groups (or even of finitely generated groups) never fall into
certain categories of groups which are very close to being free, thus they 
need to be studied directly rather than relying on results already obtained
for these special categories. In Section 3 we introduce our process
of directly finding presentation 2-complexes which are non positively
curved square complexes by starting with a 2-generator 1-relator presentation
of even length where one generator has zero exponent sum in this relator.
We then replace this generator by rewriting the relator in terms of conjugates
of the other generator. This introduces lots of length 4 relators but the
original relator will now almost certainly be much longer. However as it still
has even length, we can reduce it by two each time on substitution of new
generators that are set equal to length 3 subwords of the long
relator. This
process always terminates with a presentation where all relators have length
4 and Proposition 3.1 gives a general criterion that ensures the resulting
presentation 2-complex is non positively curved.

The rest of Section 3 then gives a particular case of this construction where
the original 2-generator 1-relator presentation can be chosen so that the
group is a strictly ascending HNN extension of a finite rank free group
(by use of K.\,S.\,Brown's algorithm for 1-relator groups) and satisfies
the $C'(1/6)$ small cancellation condition (thus is word hyperbolic), giving
us a specific strictly ascending HNN extension of a finite rank free group
that is linear over $\Z$. (Although we point out that it
is only known to be linear over $\Z$ by invoking the work of Wise and Agol:
there is no idea as to what a faithful representation in $GL(n,\Z)$ could
look like.)

In Section 4 we look at other approaches. It was shown in \cite{wssc} that
any group with a finite presentation satisfying the $C'(1/6)$ small
cancellation condition acts properly and cocompactly on a CAT(0) cube
complex. This result allows us to observe that generic strictly
ascending HNN extensions of finite rank free groups, where we take a
random endomorphism of length $n$, are not just word hyperbolic but are
linear over $\Z$, which was surely the opposite of the expected picture
before the work of Wise and Agol.

We also give an example of a straightforward injective
endomorphism of $F_2$ where the resulting strictly ascending HNN extension
has a presentation 2-complex which is a non positively curved squared complex
and where we can show this group is word hyperbolic directly from the 
endomorphism without recourse to small cancellation results. To do this,
we recall that if $\theta:F_n\rightarrow F_n$ is an injective endomorphism
then there is a necessary condition for word hyperbolicity of the resulting
strictly ascending HNN extension which is: we cannot have a non identity
element $w$ of $F_n$ and $i,j>0$ such that $\theta^i(w)$ 
is conjugate in $F_n$
to $w^j$ (which would introduce Baumslag Solitar subgroups). Although it is
still open as to whether this is a sufficient condition in general, it was
shown in \cite{kap} that it is sufficient
when $\theta$ is an immersion, meaning there
is no cancellation in the image of $\theta$. Our example is an immersion
and although we still need to show that the above condition holds, we 
establish this directly in Theorem 4.2, thus giving us a word hyperbolic 
strictly ascending HNN extension of a finite rank free group with an 
especially short description and which is linear over $\Z$.

A common test case for such groups is called the Sapir group and given by
the endomorphism $a$ maps to $ab$ and $b$ to $ba$. As this is an
immmersion, we also show directly in Theorem 4.1 that the condition in
\cite{kap} for word hyperbolicity applies to this group too. This gives us
the first written proof that the Sapir group is word hyperbolic, as the
others were all word of mouth.

In the last section we note other properties that are obtained from
applying the Agol - Wise work. Unlike linearity over $\Z$ we may need to drop
to a finite index subgroup, but we can intersect finite index subgroups to 
ensure that all our strong properties appear in the same finite index
subgroup, as listed in Theorem 5.2. For instance we obtain word hyperbolic
strictly ascending HNN extensions of finite rank free groups where all
quasiconvex subgroups are separable, but not all finitely generated or even
finitely presented subgroups are.

\section{Free-like groups}

Three well studied classes of finitely presented groups which exhibit
similar behaviour to (non abelian) free groups are: free products 
(assumed non trivial
and not equal to $C_2*C_2$), groups of deficiency at least two, meaning they
possess a finite presentation where the number of generators minus the number
of relators is greater than or equal to 2, and (non abelian) limit groups.
To bring out this comparison we can use
properties of $L^2$-Betti numbers, though here we will only need the first
$L^2$-Betti number $\beta_1^{(2)}(G)\geq 0$. 
In particular we have the following
standard facts for $G$ a finitely generated infinite group:\\
\begin{prop}
\hfill\\
(a) $\beta_1^{(2)}(H)=[G:H]\beta_1^{(2)}(G)$ for $H$ a finite index subgroup
of $G$.\\
(b) $\beta_1^{(2)}(F_n)=n-1$.\\
(c) $d(G)-1\geq \beta_1^{(2)}(G)$ where $d(G)$ is the minimum
number of generators for $G$.\\
(d) If $G=G_1*G_2$ then $\beta_1^{(2)}(G)=1+\beta_1^{(2)}(G_1) 
+\beta_1^{(2)}(G_2)-1/|G_1|-1/|G_2|$, where we interpret $1/|G|$ as 0 if
$G$ is infinite.\\
(e) If $G$ is finitely presented then $\beta_1^{(2)}(G)\geq \mbox{def}(G)-1$.\\
(f) If $G=\langle x_1,\ldots , x_n|r\rangle$ for $r$ a non trivial relator
then $\beta_1^{(2)}(G)=n-1/k$ where $r=w^k$ and $w$ is not a proper power,
by \cite{dksln}.\\
(g) A non abelian limit group must have
strictly positive first $L^2$-Betti number.
\end{prop}
The last part comes from
\cite{pichot} which established semi-continuity of $\beta_1^{(2)}$ on the
space of marked groups and then used the result of \cite{ch} that 
(non abelian) limit groups are limits of (non abelian)
free groups in this space.

If a finitely generated group $G$ possesses a homomorphism onto
$\Z$ (equivalently if the ordinary first Betti number $\beta_1(G)>0$)
then we can ask whether its kernel is finitely generated. If 
$b=\beta_1(G)=1$ then this is a straight yes/no question because there are
only two such homomorphisms, one the negative of the other and both sharing
the same kernel,
but if $b\geq 2$
then $G$ has infinitely many surjective homomorphisms $\chi$
to $\Z$ and we can ask this question for each such $\pm\chi$. The BNS
invariant $\Sigma(G)\subseteq S^{b-1}$ as introduced in \cite{bns} is a subset
of this $b-1$
dimensional sphere, proven in that paper to be an open subset, where
one views the rationally defined points of $S^{b-1}$ as corresponding to
surjective homomorphisms $\chi:G\rightarrow\Z$ and it was shown that
$\chi\in\Sigma\cap-\Sigma$
if and only if the kernel $K$ of $\chi$ is finitely generated, in
which case $G$ is the semidirect product $K\rtimes\Z$.

If $\chi\in\Sigma$ but not in $-\Sigma$ then by \cite{bns} 
Proposition 4.3 $G$ can
be written as a strictly ascending HNN extension with finitely generated
base $B$. This means that there is an injective but non surjective
homomorphism $\theta$ of $B$ so that $tBt^{-1}=\theta(B)$ for $t$ the stable
letter of the extension, whereupon $G=\langle t,B\rangle$ and 
$K=\cup_{i\in\N}t^{-i}Bt^i$. Conversely if $G$ can be written as
a strictly ascending HNN extension with finitely generated
base $B$ and $\chi$ is a homomorphism from $G$ onto $\Z$ with $B$ in the
kernel of $\chi$ then $\chi$ is in $\Sigma$ but not in $-\Sigma$.
Thus openness of $\Sigma(G)$ implies that it is
empty if and only if $G$ has no expression as an ascending HNN extension with
finitely generated base $B$. Here ascending includes both strictly ascending
and being a semidirect product $K\rtimes\Z$.

Returning to our three classes of groups, we have that a free product
surjects to $\Z$ if either of the factors do, whereas groups of deficiency
at least 2 and non abelian limit groups have first Betti number at least 2
so we can ask about their BNS invariants. In fact it is well known that no
group in any of these three classes can have a homomorphism to $\Z$ with
finitely generated kernel. This follows immediately from the fact that
these
groups have no finitely generated infinite normal subgroups of infinite
index. 
In the case of free products this is Lemma 11.2 in \cite{hmp}, but
for all three classes we can use \cite{gbsurlaL2} by Gaboriau which shows that
a finitely generated group $G$ with $\beta_1^{(2)}(G)>0$ has this property.
Alternatively this lack of finitely generated
normal subgroups was shown in 
\cite{brhw} to hold for limit groups, by establishing 
a result which generalises
that for free products: if we have a
finitely generated normal subgroup in an amalgamated free product or HNN
extension and $A$ is the amalgamated subgroup, or one of the associated
subgroups, then either $N\subseteq A$ or $AN$ has finite index in $G$.

However we do not know of an argument in the literature
showing that these
groups cannot be expressed as a strictly ascending HNN extension with finitely
generated base, so we present a straightforward proof here using rank
gradient, defined in \cite{lac} for finitely generated groups $G$ as
\[RG(G)=\inf_{H\leq_f G}\frac{d(H)-1}{[G:H]}\] 
(in fact only finite index normal subgroups were used there but it gives
the same value). It was also shown there that if $G$ is an ascending
HNN extension with finitely generated base $B$ and stable letter $t$ then
$RG(G)=0$ because we can take the finite cyclic covers
$\langle t^k,B\rangle$ which have index $k$ in $G$ and are generated by
$d(B)+1$ elements. We need to relate this quantity to 
$\beta_1^{(2)}(G)$ which
can be done using cost and/or L\"uck's approximation theorem, but the 
following is a quick argument which does not require residual finiteness.
\begin{thm}
For $G$ any finitely generated infinite
group, we have $RG(G)\geq\beta_1^{(2)}(G)$.
\end{thm}
\begin{proof}
As $d(G)-1\geq\beta_1^{(2)}(G)$ by (c), we use (a) to obtain for any
$H\leq_f G$
\[\frac{d(H)-1}{[G:H]}\geq \frac{\beta_1^{(2)}(H)}{[G:H]}=\beta_1^{(2)}(G)
.\]
\end{proof}
\begin{co}
No group $G$ which is a non trivial free product or has a presentation of
deficiency at least two or is a non abelian limit group can be expressed as
an ascending HNN extension with finitely generated base, thus the
BNS invariant $\Sigma(G)$ is empty.
\end{co}
\begin{proof}
This would require $G$ to be finitely generated with $RG(G)=0$ but 
$RG(G)\geq\beta_1^{(2)}(G)>0$ by Theorem 2.2.
\end{proof}

\section{A linear strictly ascending HNN extension}

Remarkable recent work has shown the extraordinarily strong properties
possessed by all fundamental groups of closed hyperbolic 3-manifolds $M$.
Indeed by \cite{km}, \cite{bws}, \cite{ag} and \cite{wslng} we have that
$\pi_1(M)$ is virtually special, meaning that it has a finite index
subgroup which embeds in a Right Angled Artin Group (a RAAG).
Consequently any group theoretic property which holds for all RAAGs and
is preserved by subgroups and finite index supergroups will hold for a
virtually special group. In particular $\pi_1(M)$ is linear over $\Z$,
which is a very strong property
and was not previously suspected to hold for hyperbolic
3-manifold groups. However in the course of showing that these groups
$\pi_1(M)$ are virtually special, one requires the Kahn-Markovic result
in \cite{km} that there are many surface subgroups in $\pi_1(M)$, giving
enough codimension 1 subgroups for \cite{bws} to show that $\pi_1(M)$ acts
properly and cocompactly on a CAT(0) cube complex. This then allows
application of the Agol result which tells us that any group $G$ which is
word hyperbolic and acts properly and cocompactly on a CAT(0) cube
complex has a virtual quasiconvex hierarchy. This conclusion feeds into
the machinery of \cite{wslng} to show that $G$ is virtually special. Also in
\cite{wslng} hyperbolic limit groups and 1-relator groups with torsion were
shown to have such hierarchies. Thus one might wonder what can occur for
an arbitrary 1-relator group. If this group $G$ is not word hyperbolic
(such as the Baumslag-Solitar groups) then we do not have such strong
conclusions holding in general. However if it is then we might try to
show that the hypotheses for Agol's result hold, though we do not have a
direct equivalent of the Kahn-Markovic work here. On the other hand we might
try to show that Agol's result cannot apply ($G$ ``cannot be cubulated'')
by establishing that one of the conclusions does not hold for $G$, and
linearity over $\Z$ seems strong enough to find counterexamples.

A well known necessary
condition for linearity over $\Z$ is that every soluble subgroup
must be polycyclic. In particular this rules out the soluble
Baumslag-Solitar groups $BS(1,k)=\langle x,y|yxy^{-1}=x^k\rangle$ for $k\neq 0,
\pm 1$. But these are all strictly ascending HNN extensions with base
a free group, where this free group has rank 1. Thus it might be asked
whether a strictly ascending HNN extension $E$ of a finitely generated
group $B$ can ever be linear over $\Z$ (the free group $F_2$
can be expressed as a strictly
ascending HNN extension with infinitely generated base, so finitely
generated is needed here). Note that 3-manifold groups can never be
strictly ascending HNN extensions with finitely generated base, because
their BNS invariants are symmetric. Alternatively we can specialise to when
this finitely generated base is a free group and here we can ask whether
$E$ is linear over $\C$ or even residually finite. 
The latter question was solved positively
in \cite{brsp} using algebraic geometry. As for the
former, an example which is not linear was given in \cite{ds} whereas
the groups $BS(1,k)$ embed in $\slt$. However these two
examples contain (or are) a soluble Baumslag-Solitar group, so cannot be
word hyperbolic or linear over $\Z$. Indeed Problem 6 of \cite{ds} asks
whether there are word hyperbolic non linear ascending HNN extensions
of finitely generated free groups and conjectures that if the extension
is strictly ascending then most of these groups are non linear. However
until now we are unaware of any word hyperbolic strictly ascending
HNN extensions of a finitely generated group which have been shown either
to be linear or to be non linear. Moreover we did not know of a 
specific strictly
ascending HNN extension of any finitely generated group which has shown
to be linear over $\Z$. (By Proposition 2.1, 1-relator groups
with torsion and limit groups have positive first $L^2$-Betti number 
so have empty BNS invariant by Corollary 2.3.)

We proceed by looking for particular finite
presentations that give rise to
non positively curved square complexes, which will allow us to use the
Agol result if the group given by this presentation is also
known to be word hyperbolic. Therefore suppose we have
a finite presentation of some
group $G$ where all relators are cyclically reduced and of length exactly
4. We can then form a cube complex, or more specifically a square complex,
by taking one vertex, a loop for each generator and then a square for each 
relator,
glued to form the usual presentation 2-complex. It is a standard
exercise in the theory of cube complexes to show that a square complex
is non positively curved (under the metric realising each square as
the unit square in $\R^2$ and then gluing) if and only if the length of the
shortest loop in each graph formed by taking the link of a vertex is at
least 4. If so then $G$ is the fundamental group of this non positively
curved compact 2-complex and so will act properly and cocompactly on its
universal cover, which is simply connected and inherits the non positively
curved condition, thus is a CAT(0) cube complex. Hence on showing word
hyperbolicity for $G$, we directly know that Agol's theorem applies
without having to go through the process of ``cubulating'' $G$.

As for the non positively curved condition, we only have one vertex
to deal with. Using the reformulation in terms of the link at this
vertex, we see that this is equivalent to our presentation
satisfying both of the following 
two properties:\\
\hfill\\
(1) (No length two loops)
If $L$ is the list of all cyclic permutations of the relators and their
inverses and $xy$ is a length 2 subword of something in $L$ then $xy$ appears
only once in all of $L$.\\
\hfill\\
(2) (No length three loops)
If $xy$ is a length two subword appearing somewhere in $L$, as is
$y^{-1}z$ then we cannot have $xz$ as a subword of something in $L$, other
than if $z=x^{-1}$ and these two length 4 words in $L$ are self inverse.\\
\hfill\\
These conditions are similar to, and for our presentations with all
relators of length 4 are equivalent to,  
the $C(4)$ and $T(4)$ conditions in combinatorial
group theory. A piece is a subword of an element $l$ in $L$ which cancels
completely when $l$ is multiplied by an element in $L$ which is not $l^{-1}$.
Thus (1) in our case is asserting that all pieces are of length 1.
The $C(n)$ condition in general states that every relator is a product of at 
least $n$ pieces and, though we do not provide a definition of the $T(n)$
condition as it will not be needed, we note here that the $T(4)$ is
equivalent to (2) for our presentations.
Such presentations have been studied before, for instance 
they hold for a particular
presentation (the Dehn presentation) of the fundamental group of a prime
alternating knot. However a knot group contains $\Z\times\Z$ and so cannot
be word hyperbolic. 

Here we try a more direct approach: we want to start with
a 1-relator presentation and end up with a presentation consisting
only of relators of length 4. If this relator $r=x_1\ldots x_{2n}$ is of
even length, where each $x_i$ is a generator or its inverse, then there is an
obvious approach to try: we can introduce a new generator, $u_1$ say, and set
$u_1=x_nx_{n+1}x_{n+2}$, where we assume for now that
this length 4 relator satisfies our
two conditions above. On substituting $u_1$ into $r$, we reduce the word
length of $r$ by 2. We then set the next new generator $u_2$ equal to
$x_{n-1}u_1x_{n+3}$, and these two relators have a good chance of being
valid because $u_2$ is new and $u_1$ has only appeared twice. We can
continue this process, introducing a new generator and reducing the
length of $r$ by two each time until we are left with a presentation
where all relators are of length 4, with the final relator being
$x_1x_2u_{n-2}x_{2n}$, which we might hope to pass (1) and (2).

However this approach is doomed to fail the non positively curved
condition if we have only two generators $a$ and $t$, because of our need
to pick the starting and final length 4 relator carefully. For instance
$u_1=a^{-1}ta$ and $u_{n-2}=ata^{-1}$ would pass on their own.
But this uses up all two letter subwords involving both $a^{\pm 1}$ and
$t^{\pm 1}$, yet if $n>4$ then we also have $u_2=x_{n-1}u_1x_{n+3}$ for
$x_{n-1},x_{n+3}\in\{a^{\pm 1},t^{\pm 1}\}$, thus meaning that
$au_1(x_{n-1}u_1)^{-1}=ax_{n-1}^{-1}$ has already appeared and we are now
stuck, unless $x_{n-1}=a^{-1}$ as well. Continuing this argument, we see that
$r=t^{-1}a^{1-n}ta^{n-1}$ is the only relator that passes, which will
not give rise to a word hyperbolic group. Hence
we call on a standard technique in the theory of 1-relator groups: we alter
the presentation so that one of the generators, here $t$, has zero
exponent sum in the relator. We then set $a_n=t^nat^{-n}$ for $n\in\Z$, where
$a$ is the other generator and rewrite the relator so that our presentation
is of the form
\[\langle t,a_k,\ldots ,a_0,\ldots ,a_l|w(a_k,\ldots ,a_l),
a_i=ta_{i-1}t^{-1}
\mbox{ for } k+1\leq i\leq l\rangle.\]
This introduces lots of extra generators and length 4 relators. However we now
have to go through the above process of bringing in yet more
generators to shrink $w$, which we will now suppose has even length. We will
further suppose that $k=0$ (this can be achieved by replacing $a$
with a suitable conjugate).

We can now give a general result as to when
this process produces a non positively curved
cube complex.
\begin{prop}
Suppose that $\langle a,t|r(a,t)\rangle$ is a 2-generator 1-relator
presentation of a group $G$ with $t$ having exponent sum 0 in $r$ and
set $a_i=ta_{i-1}t^{-1}$ for $i\in\N$. Suppose further that $r$ has even
length $2n$, can be
written solely in terms of $a_0,\ldots , a_l$  for some $l\geq 2$, and is
of the form
\[a_l^{-1}a_0^{-1}x_3\ldots x_{n-1}a_0^{-1}a_la_0x_{n+3}\ldots x_{2n-1}a_0\]
where all the unspecified letters $x_i$ can be picked arbitrarily from
$\{a_1^{\pm 1},\ldots ,a_{l-1}^{\pm 1}\}$ as long as $r$ is a reduced word.
Let us introduce new generators $u_1,\ldots ,u_{n-2}$ by
$u_1=a_0^{-1}a_la_0,u_2=x_{n-1}u_1x_{n+3}, \ldots,
u_{n-2}=x_3u_{n-3}x_{2n-1}$ and substitute them successively in $r$ until
we are left with $a_l^{-1}a_0^{-1}u_{n-2}a_0$. This produces a
presentation for $G$ consisting of the $l+n$ generators
$t,a_0,\ldots ,a_l,u_1,\ldots ,u_{n-2}$ and the $l+n-1$ length 4 relators
$ta_{i-1}t^{-1}a_i^{-1}$ for $1\leq i\leq l$ and
\[a_0^{-1}a_la_0u_1^{-1}
,\, x_{n-1}u_1x_{n+3}u_2^{-1}, \ldots,
x_3u_{n-3}x_{2n-1}u_{n-2}^{-1},\, a_l^{-1}
a_0^{-1}u_{n-2}a_0\] whose
presentation 2-complex is a non positively curved squared complex.
\end{prop} 
\begin{proof}
We first check condition (1) and consider when there is a repetition of a
generator or its inverse in these relators. There is no problem with the
$u_i$s; for instance $u_2$ only appears in the inverse of the second
relator as $u_2x_{n+3}^{-1}u_1^{-1}x_{n-1}^{-1}$ and in the third as
$x_{n-2}u_2x_{n+4}u_3^{-1}$ so we would only have a repeat of a 2 letter
subword involving $u_2$ if $x_{n+3}^{-1}=x_{n+4}$ or $x_{n-1}^{-1}=x_{n-2}$
but then $r$ is not reduced.

Although there are lots of repetitions of $t^{\pm 1}$, the preceding and
following letter is always different. As for the letters $a_i^{\pm 1}$,
the only place any two of these appear next to each other is in the two
relators $a_0^{-1}a_la_0u_1^{-1}$ and $a_0^{-1}u_{n-2}a_0a_l^{-1}$
occurring at the start and end of the substitution process applied to $r$.
But these two relators were initially chosen specifically to avoid
repetition. Otherwise any appearance of any letter $a_i^{\pm 1}$ is next
to a letter of the form $t^{\pm 1}$ or $u_j^{\pm 1}$ and we have just seen
that there are no two letter subwords which are repeats of these.

As for satisfying condition (2), we consider which letters can be cancelled
through a repeat. On considering
$ta_{i-1}^{\pm 1}t^{-1}a_i^{\mp 1}$ and $ta_{j-1}^{\pm 1}t^{-1}a_j^{\mp 1}$
with subwords having product $a_i^{\pm 1}a_j^{\pm 1}$, it soon becomes
clear that the only two letter words of the form $a_i^{\pm 1}a_j^{\pm 1}$
that are allowed in our relators are when $i=0$ and $j=l$ or vice versa,
thus explaining why the choice for the start and end substitutions has to be
of this form. As for a pair of two letter subwords involving cancellation
in $u_i^{\pm 1}$, the resulting product $a_i^{\pm 1}a_j^{\pm 1}$ must have
at least one of $i$ and $j$ between 1 and $l-1$ because the letters $x_k$
were chosen from this range, so this product does not occur as a
subword in our relators.

As for cancellation of $a_i^{\pm 1}$, this only leaves $t^{\pm 2}$ if both
relators are taken from the first $l$ relators (the ones involving $t$).
If one relator is taken from this set and one from the latter $n-1$
relators involving the $u_i$s then the resulting product would contain
$t^{\pm 1}$ and some $u_i^{\pm 1}$, which never occur together (or we could
have $ta_0\cdot a_0^{-1}a_l=ta_l$, or similarly $ta_l^{-1}$ or 
$t^{-1}a_0^{\pm 1}$ which do not occur because $a_l$ and $a_0$ are the top
and bottom of the relators $a_i$). Thus we are only left with cancellation
of some $a_i^{\pm 1}$ amongst these last $n-1$ relators. This could only
leave a product of the form $u_i^{\pm 1}u_j^{\pm 1}$ but no $u_i$s appear
next to each other, with the exception again of the two end relators
involving $a_0^{\pm 1}, a_l^{\pm 1}$ and here it can be checked directly
that no repeated length two subword occurs.
\end{proof}   

Thus we have a presentation of our original group with all relators of
length 4 and such that the universal cover of the presentation 2-complex is
a CAT(0) complex on which $G$ acts properly and cocompactly. We now 
go back to the original
2-generator 1-relator presentation $\langle a,t|r\rangle$ and make it
satisfy some form of small cancellation condition that is known
to imply word hyperbolicity.
For this we can use the $C'(1/6)$ small cancellation condition, which states
that any piece of the
cyclically reduced conjugates of $r$ and $r^{-1}$ has word length strictly
less than a 1/6th times that of $r$.
For definiteness we
let $l$ be even and set $n=l+2$, choosing the remaining letters of
$r$ in a regular fashion but not so regular as to make $r$ a word with
large pieces. Moreover we should ensure that $\beta_1(G)=2$ in order to
have many homomorphisms onto $\Z$, which requires
our relator to have zero exponent sum in $a_0,\ldots , a_l$.
Therefore we let 
\begin{eqnarray*}
r_l&=&a_l^{-1}a_0^{-1}a_{l-1}a_1a_{l-2}a_2\ldots a_{l/2}a_0^{-1}\cdot\\
   & &a_la_0a_{l-1}^{-1}a_1^{-1}a_{l-2}^{-1}a_2^{-1}\ldots a_{l/2}^{-1}a_0.
\end{eqnarray*}

We find (by hand and then by writing a simple program to check) that
for $l=8$ the relator has length 116 and the largest piece length 17 so we
have a word hyperbolic group. This gives us the relator
\[
t^8a^{-1}t^{-8}a^{-1}t^7at^{-6}at^5at^{-4}at^3at^{-2}atat^{-4}a^{-1}\cdot\]
\[
t^8at^{-8}at^7a^{-1}t^{-6}a^{-1}t^5a^{-1}t^{-4}a^{-1}
t^3a^{-1}t^{-2}a^{-1}ta^{-1}t^{-4}a, \]
from which we obtain:
\begin{thm}
The 2-generator 1-relator group with relator as above\\
(a) is a strictly ascending HNN extension of a finite rank free group,\\
(b) is word hyperbolic,\\
(c) and acts properly and cocompactly on a CAT(0) cube complex.\\ 
Thus by \cite{ag} and \cite{wslng} this group is linear over $\Z$.
\end{thm}
\begin{proof}
We now just need to establish (a). For this we use K.S.Brown's algorithm in
\cite{ksb} where we draw out the relator and take its convex hull
$C$. Noting that all vertices of $C$ have been passed through once we
conclude that the homomorphisms $\chi(a)=\chi(t)=-1$ and $\chi(a)=-\chi(t)=1$
express $G$ as a strictly ascending HNN extension of a finite rank free
group. In fact we can also conclude all other homomorphisms except
$\chi(a)=\pm 1,\chi(t)=0$ and $\chi(a)=0,\chi(t)=\pm 1$ express $G$ as a
(finite rank free)-by-$\Z$ group. Thus such groups can also be linear over
$\Z$, though this does not seem such a big surprise compared to the case 
of strictly ascending HNN extensions.   
\end{proof}

Thus we have the following corollaries:
\begin{co}
There exists a strictly ascending HNN extension of a finitely generated
(or even finite rank free) group which is linear over $\Z$.
\end{co}
\begin{co}
There exists a word hyperbolic
strictly ascending HNN extension of a finitely generated
(or even finite rank free) group which is linear over $\C$.
\end{co}

\section{Other results}

In the last section we showed that our group $G$ acted properly and
cocompactly on a CAT(0) cube complex
by explicitly exhibiting a presentation 2-complex for $G$
and then we used a small cancellation
condition to obtain word hyperbolicity. In fact further work of Wise
allows one to gain both conclusions from some small cancellation
conditions, although the CAT(0) complex constructed will not be a square
complex obtained directly from the presentation but a much higher
dimensional complex formed by considering the intersections of various
hypergraphs. Theorem 1.2 in \cite{wssc} states that if $G$ is a group
with a finite $C'(1/6)$ presentation then $G$ acts properly and
cocompactly on a CAT(0) cube complex, and moreover $G$ will be word
hyperbolic. (Thus the group in \cite{lnw} Proposition 3 also satisfies
the same conditions as in Theorem 3.2 although the square complex implied
by Figure 3 in that paper is not CAT(0).)

This gives rise to the observation that, rather than ``most'' strictly
ascending HNN extensions of finite rank free groups being non linear as
conjectured in \cite{ds}, we have that ``most'' of these groups are linear 
over $\Z$. Indeed consider the model of \cite{cw} Definition 7.4 where we
fix a free group $F_k$ and a free generating set $x_1,\ldots ,x_k$, then
regard a random endomorphism of $F_k$ of length $n$ to be given by picking
$k$ random reduced words $w_1,\ldots ,w_k$ of length $n$ and sending
$x_i$ to $w_i$. Here the words $w_1,\ldots ,w_k$ are picked independently
using the uniform distribution, and standard  counting arguments in the
theory of random groups (see \cite{oli} for an overview) tell us that the
presentation $\langle x_1,\ldots ,x_k|w_1,\ldots ,w_k\rangle$ satisfies
$C'(1/6)$, or indeed $C'(1/7)$, with probability tending to 1 as
$n\rightarrow\infty$. Here we actually have the presentation
\[\langle t,x_1,\ldots ,x_k|tx_1^{-1}t^{-1}w_1,\ldots ,
tx_k^{-1}t^{-1}w_k\rangle\]
but a relator in the latter presentation is two letters longer than in the
former, whereas a piece can only be at most two letters longer.
Consequently as $n\rightarrow\infty$ we have that the latter presentation
is $C'(1/6)$, thus is word hyperbolic and linear over $\Z$, with probability
tending to 1.

The purpose of \cite{cw} was to give conditions under which a
strictly ascending HNN extension of a finite rank free group contains a
surface group. A good test case is what has come to be known as the
Sapir group, defined using the $F_2$ endomorphism $a\mapsto ab,
b\mapsto ba$. This was shown to contain a surface group in \cite{cw} by an
implementation of an algorithm described in that paper involving linear
programming. We did not manage to show that the Sapir group
acts properly and cocompactly on a CAT(0) cube complex and so cannot
confirm that it is linear over $\Z$. However we do have an example of a
very similar looking group where this holds, by finding a suitable
non positively curved square complex directly from its presentation. We
then would like to establish word hyperbolicity without recourse to theorems
on small cancellation presentations or other such results. We can do this
by using Theorem A in \cite{kap} which considered the word hyperbolicity of
strictly ascending HNN extensions of a finite rank free group $F(X)$ 
and came up
with a result when the endomorphism $\theta$
is an immersion, which is defined to mean
that for all $x,y\in X\cup X^{-1}$ with $xy\neq e$, the word 
$\theta(x)\theta(y)$ admits no cancellation. This implies that if
$w=x_1\ldots x_n$ is written as a reduced word on $X$ for $x_i\in
X\cup X^{-1}$ then no cancellation can occur on the right hand side of
$\theta(w)=\theta(x_1)\ldots \theta(x_n)$. Thus the word length of
$\theta(w)$ is the sum of that for the $\theta(x_i)$. This theorem states that
if there is no periodic conjugacy class, meaning that there is no
$w\in F(X)-\{e\}$ and  $i,j>0$ such that $\theta^i(w)$ is conjugate to
$w^j$ in $F(X)$, then the strictly ascending HNN extension $F(X)*_\theta$
is word hyperbolic. For instance the Sapir group endomorphism is an
immersion. We note that no one has been inclined to write down a proof
that this group is word hyperbolic (for instance in \cite{ds} it was stated
that this was by \cite{kap} but in \cite{spgp} this became ``it is known
(Minasyan) that this group is hyperbolic'' without reference). Thus it
seems reasonable to provide a proof here which does use I.\,Kapovich's
Theorem.
\begin{thm}
The Sapir group
\[F(a,b)*_\theta=\langle t,a,b|tat^{-1}=ab, tbt^{-1}=ba\rangle\]
is word hyperbolic.
\end{thm}
\begin{proof}
By the above mentioned Theorem of Kapovich, we need to show that if
$\theta^i(w)$ is conjugate to $w^j$ for $i,j>0$ then $w=e$.
Now $\theta$ 
being an immersion and each generator mapping to a word of
length two implies that the word length of $\theta^i(w)$ is $2^i$ times
that of $w$.
But we can assume by conjugating that $w$ is cyclically reduced, in which
case $\theta^i(w)$ is also cyclically reduced, as will be $w^j$ whose
word length is $j$ times that of $w$. Now
if $\theta^i(w)$ and $w^j$ are both cyclically reduced and are conjugate
in $F_2$ then they must have the same length and be cyclic permutations
of each other, thus $j=2^i$. 

Next suppose the word length of $w$ is even and set $w=uv$ for $u,v$ both
half the length of $w$. We have $\theta^i(u)\theta^i(v)$ is a cyclic
permutation of $w^{2^i}$ and so is equal to $(sr)^{2^i}$
without cancellation in this expression, where $w=rs$ (but $r$ and $s$ are
not necessarily of equal length). Thus 
$\theta^i(u)=(sr)^{2^{i-1}}=\theta^i(v)$. 
As an immersion must be injective, we obtain $u=v$ and $\theta^i(u)$
is a cyclic permutation of $w^{2^{i-1}}=u^{2^i}$.
Thus we can
continue to cut $w$ in half at each stage until it has odd length, whilst
still preserving the fact that $\theta^i(w)$ is a cyclic conjugate of
$w^{2^i}$. 
When we reach this point, we
consider the number of appearances $a^{\pm 1}$ and $b^{\pm 1}$ in $w$, which
cannot be equal, thus nor can it be equal for $w^{2^i}$ or for any cyclic
conjugate thereof.
But as $\theta(a)$ and
$\theta(b)$ both have equal exponent sums, this is true for any element in
the image of $w$ which is a contradiction.
\end{proof}

We now give our closely related example.
\begin{thm}
The group
\[F(a,b)*_\phi=\langle t,a,b|tat^{-1}=ab^{-1}a^2b, 
tbt^{-1}=ba^{-1}b^2a\rangle\]
acts properly and cocompactly on a CAT(0) square complex and is word 
hyperbolic, thus is linear over $\Z$.
\end{thm}
\begin{proof}
We take the universal cover of
the square complex below in Figure 1 which satisfies conditions
(1) and (2) in Section 3 (where we have placed some arbitrary extra letters
for the new generators that are introduced).
Moreover the endomorphism $\phi$ is an
immersion, so we can provide a similar proof to that
of Theorem 4.1 which will show our group is word hyperbolic. 
We can assume we have 
$\phi^i(w)$ being cyclically conjugate to $w^{5^i}$ by the same
argument as before because here words are mapped to five times their
length under $\phi$. 
The exponent sums of $a$ and $b$ in $\phi(a)$ are
3 then 0, and reversed for $\phi(b)$. Thus if we set $\lambda$ for the
exponent sum of $a$ in $w$ then the exponent sum of $a$ in $\phi(w)$
must be $3\lambda$ and $3^i\lambda$ in $\phi^i(w)$
but it must be $5^i\lambda$ for $w^{5^i}$. However these elements are
conjugate so the exponent sums must be equal, forcing $w$ to have
zero exponent sum in $a$.

We can assume that $w$ contains the letter $a$ or $b$, as an appearance
of a negative letter in $w$ produces positive letters in $\phi(w)$ and
positive letters will continue to be present under iteration of $w$ until
we obtain $\phi^i(w)$. This is a cyclic permutation of a power of $w$
and thus each letter of $\phi^i(w)$ occurs somewhere in $w$. Moreover
as $\phi$ commutes with the 
involution $\iota(a)=b$, $\iota(b)=a$, we can say that $w$ contains the
letter $a$ without loss of generality. 

Now consider the semi-infinite word
\[s=ab^{-1}a^2ba^{-1}b^{-2}ab^{-1}ab^{-1}a^2bab^{-1}a^2b\ldots\]
formed by repeated iteration of $\phi$ starting from $a$. As $a$ appears
somewhere in $w$ and there is no cancellation when applying $\phi$, the
word $\phi^i(a)$ (which is the first $5^i$ letters of $s$) appears in
$\phi^i(w)$ and so appears as a cyclic subword in the cyclic word given by
$w^{5^i}$. If the length of $\phi^i(a)$ is at most that of $w$ then
$\phi^i(a)$ is then also a cyclic subword of the cyclic word $w$, otherwise
$\phi^i(a)$ will ``wrap around $w$''. In the former case we must have that
$\phi^i(a)$ being a cyclic subword of $w$ implies that 
$\phi^{2i}(a)=\phi^i(\phi^i(a))$ is a cyclic subword of $\phi^i(w)$, and thus
of $w^{5^i}$. We can now iterate this process until we cover $w$, whereupon
we conclude that there is an initial prefix of our semiinfinite word $s$
which is a cyclic permutation of $w$ 

We finish by showing that no initial prefix of $s$ has zero exponent sum
in $a$, to obtain a contradiction. Otherwise
let us take the shortest initial
prefix $x$ of $s$ which has 0 exponent sum in $a$, with length
$5q+r$ for $r=0,1,2,3,4$. If $r=0$ then we can take the first fifth $v$
of $x$ and we have $\phi(v)=x$, thus this has zero exponent sum in $a$ too.
Otherwise we knock off the last $r$ letters of $x$, which also results in
an element with length divisible by 5. Moreover this element has exponent
sum in $a$ at most 2, by considering the presence of $a^{\pm 1}$ in
$\phi(a^{\pm 1})$ and $\phi(b^{\pm 1})$. This gives us a contradiction
as before, as either the exponent sum of this element is not divisible
by 3 or it is at most 0.
\end{proof}
\begin{figure}[h]
\begin{center}
\includegraphics[angle=-90,width=12cm]{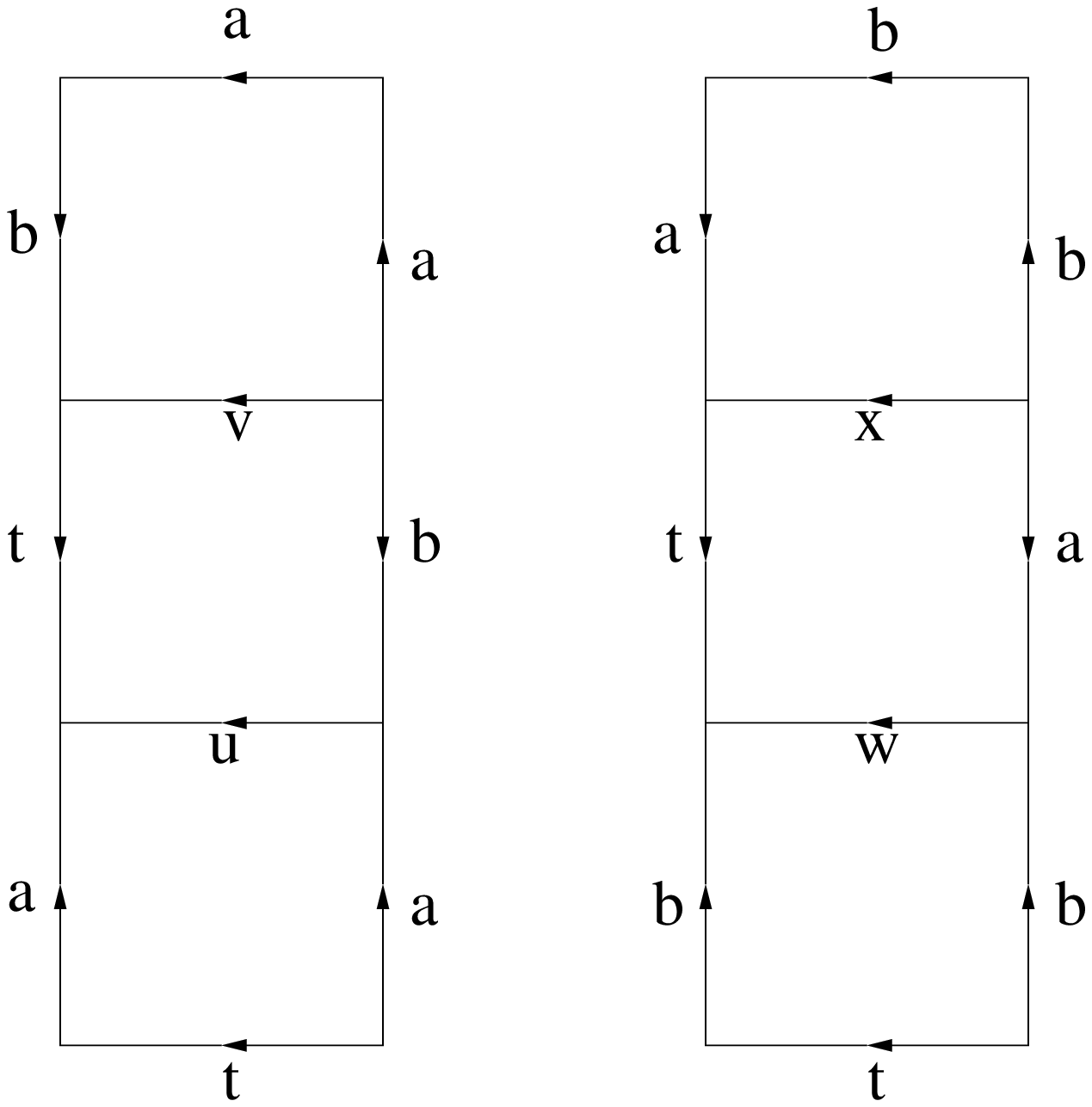}
\caption{}
\end{center}
\end{figure} 

\section{Further Properties}

Having used the Wise - Agol machinery to conclude that our strictly
ascending HNN extension $G$ from the last section is linear over $\Z$,
it is now worth seeing what other properties can be obtained by this
approach and whether they were already known. Largeness (having a
finite index subgroup surjecting onto the free group $F_2$) is also
implied, because a subgroup of a RAAG either has this property or is
free abelian by \cite{ashi}. Now some strictly ascending HNN extensions
of finite rank free groups were shown to be large in \cite{bdf1} though none
of these are word hyperbolic (in fact conjecturally they are precisely
the non word hyperbolic such groups). But \cite{bcmp} describes a
computer program that can determine largeness given a finite presentation
as input and this program instantly says yes on input of our group $G$
in Theorem 4.2.

However \cite{ashi} actually obtains a stronger property: any subgroup of a
RAAG either surjects to $F_2$ or is free abelian. Here the subgroup does not
need to be of finite index, nor even finitely generated. Thus on taking
a finite index subgroup $H$ of $G$ which is special, we have that $H$ is
also a strictly ascending HNN extension by \cite{bdf1} Proposition 4.3 (iii),
with the extreme property that every non cyclic subgroup of $H$
surjects to $F_2$
(as $H$ is also word hyperbolic). Note that $G$ itself does not surject to
$F_2$ because the abelianisation of $G$ has rank 1 over $\Z$.

Another strong property gained by this approach is being LERF, meaning that
every finitely generated subgroup is separable (an intersection of finite
index subgroups). This is because \cite{wslng} shows that
if a word hyperbolic group is virtually special then all quasiconvex
subgroups are separable. Now it was already known that a non quasiconvex
finitely generated subgroup of a closed hyperbolic 3-manifold group is
virtually a fibre in a fibration, thus putting these two facts together implies
LERF for closed hyperbolic 3-manifold groups. Yet strictly ascending HNN
extensions with finitely generated base $B$ are the archetypal examples of
non LERF groups because we cannot separate $B$ and $tBt^{-1}$.
Therefore we have
\begin{prop}
The groups $G$ in Theorems 3.2 and 4.2
are word hyperbolic groups where all
quasiconvex subgroups are separable, but not all finitely presented
subgroups are.
\end{prop}

Finally we have various algebraic properties holding for the special
finite index subgroup $H$ which follow because they are true for RAAG groups
and are preserved by subgroups. These include being
residually torsion free nilpotent by \cite{hsws} and hence being
biorderable, as well as being virtually RFRS 
(residually finite rational soluble) by
\cite{agr}. This last can be replaced by RFRS simply by dropping down
to the relevant finite index subgroup $L$ of $H$
as all other properties will be preserved. 
Note though that $H$ and $L$ will not be
residually free: if so then a result of B.\,Baumslag implies that they are 
fully residually free (as otherwise they must contain $F_2\times\Z$ but
they do not even contain $\Z\times\Z$). However this is ruled out by 
Proposition 2.1 (g) and Corollary 2.3.

Nevertheless we have a
group with remarkable properties which are summarised here:
\begin{thm} There exists a strictly ascending HNN extension $L$ of a finite
rank free group which is word hyperbolic and such that\\
(i) $L$ is linear over $\Z$.\\
(ii) All quasiconvex subgroups of $L$ are separable but not all finitely
presented subgroups.\\
(iii) $L$ is biorderable.\\
(iv) $L$ is residually torsion free nilpotent.\\
(v) $L$ is RFRS.\\
(vi) The only subgroups of $L$ not possessing a surjection to $F_2$ are the
infinite cyclic subgroups and the trivial subgroup.
\end{thm}

There are word hyperbolic groups which do not
have a finite index subgroup with such properties, for instance
the last three conditions immediately imply the existence of a
surjection to $\Z$, so a word hyperbolic group with property
(T) has no finite
index subgroup possessing any of (iv), (v), (vi). Nevertheless, we
wonder if all word hyperbolic groups in well behaved classes with no
obvious obstruction have such a finite index subgroup. As examples of
well behaved classes with no
obvious obstruction, we might take all $F_n$ by $\Z$ groups or all
strictly ascending HNN extensions of finite rank free groups, or all
1-relator groups.


\begin{thebibliography}{99}

\bibitem{agr} I.\,Agol,
{\it Criteria for virtual fibering},
J. Topol. {\bf 1} (2008), 269--284.

\bibitem{ag} I.\,Agol,
{\it The virtual Haken conjecture}, (2012)\\
\texttt{http://arxiv.org/abs/1204.2810}

\bibitem{ashi} Y.\,Antolin and A.\,Minasyan,
{\it Tits alternatives for graph products}, (2011)\\
\texttt{http://arxiv.org/abs/1111.2448}

\bibitem{bws} N.\,Bergeron and D.\,T.\,Wise,
{\it A boundary criterion for cubulation},
Amer. J. Math. {\bf 134} (2012), 843--859. 

\bibitem{bns} R.\,Bieri, W.\,D.\,Neumann and R.\,Strebel,
{\it A geometric invariant of discrete groups},
Invent. Math. {\bf 90} (1987), 451--477.

\bibitem{brsp}A.\,Borisov and M.\,Sapir,
{\it Polynomial maps over finite fields and residual finiteness of mapping
tori of group endomorphisms},
Invent. Math. {\bf 160} (2005), 341--356.

\bibitem{brhw} M.\,R.\,Bridson and J.\,Howie,
{\it Normalisers in limit groups},
Math. Ann. {\bf 337} (2007), 385--394.

\bibitem{ksb} K.\,S.\,Brown,
{\it Trees, valuations, and the Bieri-Neumann-Strebel invariant},
Invent. Math. {\bf 90} (1987), 479--504.

\bibitem{bdf1} J.\,O.\,Button,
{\it Large groups of deficiency 1},
Israel J. Math. {\bf 167} (2008), 111--140.

\bibitem{bcmp} J. O. Button, 
{\it Proving finitely presented groups are
large by computer}, 
Exp. Math. {\bf 20} (2011), 153--168.

\bibitem{cw} D.\,Calegari and A.\,Walker,
{\it Surface subgroups from linear programming}, (2012)\\
\texttt{http://arxiv.org/abs/1212.2618}

\bibitem{ch} C.\,Champetier and V.\,Guirardel,
{\it Limit groups as limits of free groups},
Israel J. Math. {\bf 146} (2005), 1--75.

\bibitem{dksln} W.\,Dicks and P.\,A.\,Linnell,
{$L^2$-Betti numbers of one relator groups},
Math. Ann. {\bf 337} (2007), 855--874.

\bibitem{ds} C.\,Dru\c tu and M.\,Sapir,
{\it Non-linear residually finite groups},
J. Algebra {\bf 284} (2005), 174--178.

\bibitem{hsws} G.\,Duchamp and D.\,Krop,
{\it The lower central series of the free partially commutative group},
Semigroup Forum {\bf 45} (1992), 385--394.

\bibitem{gbsurlaL2} D.\,Gaboriau
{\it Sur la (co-)homologie $L^2$ des actions pr\'{e}servant une mesure},
C. R. Acad. Sci. Paris S\'{e}r. I Math. {\bf 330} (2000), 365--370.

\bibitem{hmp} J.\,Hempel, {\it 3-manifolds}, Ann. of Math. Studies No. 86,
Princeton University Press, Princeton, N.\,J.\,, 1976.

\bibitem{kap} I.\,Kapovich, 
{\it Mapping tori of endomorphisms of free groups},
Comm. Algebra {\bf 28} (2000), 2895--2917.

\bibitem{km} J.\,Kahn and V.\,Markovic,
{\it Immersing almost geodesic surfaces in a closed hyperbolic
three manifold},
Ann. of Math. {\bf 175} (2012), 1127--1190. 

\bibitem{lac} M. Lackenby, {\it Expanders, rank and graph of groups},
Israel J. Math. 146 (2005), 357--370.

\bibitem{lnw} I.\,J.\,Leary, G.\,A.\,Niblo and D.\,T.\,Wise,
{\it Some free-by-cyclic groups},
Groups St. Andrews 1997 in Bath, III, 512--516, London Math. Soc. Lecture Note
Ser. 261, Cambridge Univ. Press, Cambridge, 1999.

\bibitem{oli} Y.\,Ollivier,
{\it A January 2005 invitation to random groups},
available at\\
\texttt{http://www.yann-ollivier.org/rech/publs/randomgroups.pdf}

\bibitem{pichot} M.\,Pichot,
{\it Semi-continuity of the first $l^2$-Betti number on the space of
finitely generated groups},
Comment. Math. Helv. {\bf 81} (2006), 643--652.

\bibitem{spgp} M.\,V.\,Sapir,
{\it Some group theory problems},
Int. J. Algebra Comput., {\bf 17} (2007), 1189--1214.
Ser. 261, Cambridge Univ. Press, Cambridge, 1999.

\bibitem{wssc} D.\,T.\,Wise,
{\it Cubulating small cancellation groups},
Geom. Funct. Anal. {\bf 14} (2004), 150--214.

\bibitem{wslng} D.\,T.\,Wise,
{\it From Riches to Raags: 3-Manifolds, Right-Angled Artin
Groups, and Cubical Geometry},
CBMS Regional
Conference Series in Mathematics No. 117,
American Mathematical Society, Providence, RI, 2012.


\end{thebibliography}
\end{document}